\documentclass[aps,preprint]{revtex4}
\usepackage{epsfig}
\usepackage{amsmath}
\usepackage{epsfig}

\begin{document}
\title
{Ahmed's Integral: the maiden solution}    
\author{Zafar Ahmed\\
Nuclear Physics Division, Bhabha Atomic Research Centre, Mumbai 400 085, India \\
zahmed@barc.gov.in, zai-alpha@hotmail.com}
\begin{abstract}
In 2001-2002, I happened to have proposed a new definite integral in the American Mathematical Monthly (AMM),
which later came to be known in my name (Ahmed). In the meantime, this integral
has been mentioned in mathematical encyclopedias and dictionaries and further it has also been cited and
discussed in several books and journals. In particular, a google search with the key word ``Ahmed's Integral'' throws up more than 60 listings. Here I present the maiden solution for this integral.

\end{abstract}

\maketitle

My proposal of evaluating the following integral 
\begin{equation}
\int_{0}^{1} \frac{\tan^{-1}\sqrt{2+x^2}}{(1+x^2)\sqrt{2+x^2}} dx
\end{equation}
was published  [1] in 2001, when it was thrown open to be solved  within next 6 months. Subsequently, 20 authors and two problem solving groups proposed correct solutions. The solutions to (1) by two authors, Kunt Dale (Norway) and  George L. Lamb Jr. (Arizona) have been published titled `Definitely an Integral' [2]. AMM usually prefers to publish the solutions of other solvers
than that of the proposer.

This integral is now known as Ahmed's Integral. A Google search ``Ahmed's Integral" brings more than 60 listings to view. It has been included in encyclopedias and dictionaries. Various solutions, extensions,  properties and connections of this integral  have been discussed  in a variety of ways. This integral
is very well discussed in two very interesting books on integrals  [3,4] and mentioned in another [5]. This analytically solvable integral also serves as a test model for various new methods of high precision numerical (quadratures) integrations [6].

Since in recent times this integral has evoked a considerable  attention, I propose to present the maiden solution that was sent along with the proposal [1].

Let us call the integral (1) as $I$ and use $\tan^{-1} z= \frac{\pi}{2} -\tan^{-1} \frac{1}{z}$ to split $I$ as $I=I_1-I_2$. Using the substitution $x=\tan \theta$,  we can write
\begin{equation}
I_1=\frac{\pi}{2} \int_{0}^{\pi/4} \frac{\cos \theta~d\theta}{\sqrt{2-\sin^2\theta}},
\end{equation}
which can be evaluated as $I_1=\frac{\pi^2}{12}$ by using the substitution  $\sin \theta =\sqrt{2} \sin \phi$. 
Next we use the representation
\begin{equation}
\frac{1}{a} \tan^{-1}\frac{1}{a}= \int_{0}^1 \frac{dx}{x^2+a^2} \quad (a\ne 0)
\end{equation}
to express 
\begin{equation}
I_2=\int_{0}^1 \int_0^1 \frac{dx~dy}{(1+x^2)(2+x^2+y^2)}.
\end{equation}
Further $I_2$ can be re-written as
\begin{equation}
I_2=\int_{0}^1 \int_0^1 \frac{1}{(1+y^2)} \left( \frac{1}{(1+x^2)}-\frac{1}{(2+x^2+y^2)}\right)~dx~dy
\end{equation}
\begin{small}
\begin{equation}
=\int_{0}^1 \int_0^1 \frac{dx~dy}{(1+y^2)(1+x^2)}-\int_{0}^1 \int_0^1 \frac{dx~dy}{(1+y^2)(2+x^2+y^2)}.
\end{equation}
\end{small}
Utilizing  the symmetry of the integrands and the domains for  $x$ and $y$, the second integral in (6) equals $I_2$ itself. This leads to
\begin{equation}
2I_2=\int_{0}^1 \int_0^1 \frac{dx~dy}{(1+y^2)(1+x^2)}=
\left (\int_0^1 \frac{dx}{1+x^2} \right)^2= \frac{\pi^2}{16}.
\end{equation}
Eventually, we get
\begin{equation}
I=\frac{\pi^2}{12}-\frac{\pi^2}{32}=\frac{5\pi^2}{96}. 
\end{equation}

\section*{References:}
\noindent
$[1]$ Z. Ahmed, Amer. Math. Monthly, \underline{10884},
{\bf 108} 566 (2001).\\
$[2]$ Z. Ahmed, `Definitely an Integral', Amer. Math. Monthly, \underline{10884},
 {\bf 109} 670-671 (2002). \\
$[3]$ J. M. Borwein, D.H. Bailey, and R. Girgensohn,  `Experimentation in Mathematics: Computational Paths to Discovery' (Wellesley, MA: A K Peters) pp. 17-20  (2004).\\ 
$[4]$ P. J. Nahin, `Inside Interesting Integrals' (Springer: New York) pp.190-194  (2014).\\
$[5]$ G. Boros and V.H. Moll, `Irresistible Integrals'
(Cambridge University Press: UK, USA) p.  277 (2004). \\
$[6]$ David H. Bailey, Xiaoye S. Li and Karthik Jeyabalan, `A Comparison of Three High-Precision Quadrature Programs',
Experimental Mathematics
{\bf 14}, pp. 317-329 (2005).
\end{document}